# INFERENCE FOR CENSORED QUANTILE REGRESSION MODELS IN LONGITUDINAL STUDIES

By Huixia Judy Wang[1] and Mendel Fygenson

*North Carolina State University and University of Southern California*

We develop inference procedures for longitudinal data where some of the measurements are censored by fixed constants. We consider a semi-parametric quantile regression model that makes no distributional assumptions. Our research is motivated by the lack of proper inference procedures for data from biomedical studies where measurements are censored due to a fixed quantification limit. In such studies the focus is often on testing hypotheses about treatment equality. To this end, we propose a rank score test for large sample inference on a subset of the covariates. We demonstrate the importance of accounting for both censoring and intra-subject dependency and evaluate the performance of our proposed methodology in a simulation study. We then apply the proposed inference procedures to data from an AIDS-related clinical trial. We conclude that our framework and proposed methodology is very valuable for differentiating the influences of predictors at different locations in the conditional distribution of a response variable.

**1. Introduction.** Longitudinal studies, in which repeated measurements are made on the same subject, are common in many areas of research. However, proper quantile inference procedures have not been established for longitudinal data in which some responses are left censored. This occurs, for example, when assessing the concentration of a pollutant in the environment [27], the antibody concentration in blood serum [22] or the amount of viral RNA (i.e., viral load) in individuals infected with Human Immunodeficiency Virus (HIV) [15]. In such cases, left censoring is typically due to the detection limit of the diagnostic assay.

In this paper, we consider inferences in a quantile regression setup where some of the responses are censored by fixed values and where repeated mea-

Received August 2007; revised October 2007.
[1]Supported by NSF Grant DMS-07-06963.

surements may be taken at different points across subjects. We anchor our investigation to an AIDS-related clinical trial since many approaches proposed for dealing with left censoring in longitudinal studies have been applied to such data.

Viral load is a measure of the amount of actively replicating virus and is used as a marker of disease progression among HIV-infected people. Viral load measurements are often subject to left censoring due to a lower limit of quantification. The detection limit depends upon the assay used, ranging from 500 copies/ml for the first assays available in the mid-nineties to 50 copies/ml for today's ultrasensitive assay. Despite the improvement in assay sensitivity, left censoring remains a critical issue because anti-retroviral treatments have become so effective as to lead to a steep decrease of HIV-RNA after their initiation.

Studies that measure HIV-RNA commonly incorporate repeated measurements in order to (1) control for variation among individuals and (2) monitor temporal changes in viral load during treatment. Characterization of viral dynamics in patients with different treatment regimens is essential to further development of treatments and evaluation of their efficacy (e.g., [5]).

In the medical statistical literature several methods have been proposed to handle the left censoring of HIV-RNA data. These include crude methods that use either the threshold value or some arbitrary point, such as the mid-point between zero and the cut off for the detection (e.g., [14]). These approaches usually lead to biased predictions that are systematically higher than predictions based on the true unknown values below the cut-off [8].

Other researchers considered mixed models and many applied a likelihood-based approach while assuming Gaussian distribution for both random effects and random errors; see, for example, [15, 16, 21, 32]. Chu et al. [4] considered mixture models to study the correlation between a pair of viral load measurements from each of a sample of patients assuming bivariate normal distributions. Compared to simple imputation, likelihood-based methods produce estimators that are less biased but with higher standard deviations. Even though the normality assumption eases mathematical complications, it may be unrealistic as viral load measurements are known to be highly skewed to the right, even after log transformation; see [7, 13]. Some nonlikelihood-based approaches include Sun and Wu [28], which considered a regression model with semi-parametric time-varying coefficients, and Hogan and Lee [13], which studied marginal structural quantile models with time-varying treatments. The former paper ignored the left censoring of the viral load measurements, and the latter replaced the censored values with a random generated number between zero and the detection limit.

In general, discarding censored measurements or ignoring them as such leads to biased inferences. Treating longitudinal data as independent observations can result in wrong nominal levels and/or power loss in testing hy-potheses. In what follows, we develop inference procedures within the semi-parametric framework of quantile regression. In our analysis, we examine and account for both the effects of fixed censoring in the dependent variable and the longitudinal nature of the observations. Since semi-parametric quantile regression models impose minimal assumptions on the error term, our resulting inference procedures are robust to distributional misspecifications and most appropriate for applications with extremely skewed observations. When the censored response variable is non-Gaussian, a traditional regression approach, which captures changes in the conditional mean, may not effectively detect changes in the conditional distribution. This can be critical in applications where the upper/lower quantiles of the response variable may relate differently to the covariates, leading to differing assessments of a factor's importance or a treatment's efficacy.

Powell [24, 25] pioneered inference procedures for quantile regression with fixed censoring. Bilias, Chen and Ying [2] proposed a re-sampling-based inference procedure by convexifying Powell's estimator in the resampling stage. Later, Zhao [35] discussed several median inferential methods. Ying, Jung and Wei [34] and Portnoy [23] provided quantile estimation procedures for random censoring. While other papers have been written on censored quantile regression, all existing inferential methods are developed for independent observations.

In this paper, we develop large sample inference procedures for longitudinal data. Our focus is on testing hypotheses about treatment equality and covariate significance in quantile regression models. When proposing test statistics, one may either explore the asymptotic normality of the estimated coefficients or apply likelihood ratio-based tests. However, the former requires estimating the corresponding variance–covariance matrix, which is a challenge in our semi-parametric framework because the variance–covariance matrix is a function of the unspecified densities of error terms. The latter is, in general, difficult to develop for quantile regression, and even more so in our framework because the limiting distribution takes a complicated form involving the unknown error density function. We therefore extend the rank score test proposed in [9] to our setting and study its local power theoretically and through simulations. A similar testing approach was successfully implemented in [30] in the context of conditional growth charts and in [29] for detecting differential expressions in GeneChip micoarray data.

This paper is organized as follows: In the next section we introduce notation, review various models, provide the large sample properties of the corresponding estimators, present the rank score test, and discuss the construction of confidence intervals. In Section 3 we report results from a simulation study comparing our method with two naïve methods and a bootstrap method. In Section 4, we demonstrate our method through analysis of HIV-RNA data from an AIDS clinical trial study. In Section 5, we discuss the merits of our methodology and outline future research topics. Technical proofs of the theorems and other lemmas are relegated to the Appendix.

## 2. Estimation and proposed test.

2.1. *Model setup and notation.* Longitudinal studies are typically characterized by a large number of subjects, $N$, that are each measured a relatively small number of times, $n_i$, resulting in a total of $n = \sum_{i=1}^{N} n_i$ observations. In this paper, we focus on cases where some measurements are left censored at zero. However, the proposed procedures can easily be modified to accommodate censoring from the right and/or left as long as the censoring points are fixed.

Let $y_{ij}^*$ denote the potentially left censored $j$th response of the $i$th subject, and let $y_{ij} = \max(0, y_{ij}^*)$ be its corresponding observed values. We start with the following latent regression model:

$$(2.1) \qquad y_{ij}^* = x_{ij}^T \alpha_0 + z_{ij}^T \beta + u_{ij}, \qquad i = 1, \ldots, N, \ j = 1, \ldots, n_i,$$

where $x_{ij}$ and $z_{ij}$ are the $p \times 1$ and $q \times 1$ design vectors, $\alpha_0$ and $\beta$ are $p$- and $q$-dimensional unknown parameters and $u_{ij}$ is the random error whose distribution may vary with $(x, z)$. Throughout this paper, we assume that $u_{ij}$ are independent across $i$ (subjects) but are dependent, via exchangeable correlation, within a subject. A typical example is the random intercept effect model with $u_{ij} = a_i + e_{ij}$, where $a_i$ are i.i.d. random subject effects that are independent of the i.i.d. measurement errors $e_{ij}$. We further assume that the first element of $x_{ij}$ is 1, making the first component of $\alpha_0$ an intercept.

From (2.1) and for a given $0 < \tau < 1$, we consider the following left censored quantile regression model:

$$(2.2) \quad y_{ij} = \max(0, x_{ij}^T \alpha_0 + z_{ij}^T \beta + u_{ij}), \qquad i = 1, \ldots, N, j = 1, \ldots, n_i.$$

We assume that the $\tau$th quantile of $u_{ij}$ is zero. Other than that, no distribution assumptions are made on $u$.

Since a major motivation for our study is to develop procedures for comparing HIV treatments within model (2.2), we consider testing the following hypotheses:

$$(2.3) \qquad H_0 : \beta = 0 \quad \text{versus} \quad H_n : \beta = n^{-1/2} \beta_0,$$

where $\alpha_0$ is unspecified and $\beta_0 \in \mathbb{R}^q$ is fixed. This is equivalent to comparing the null model

$$(2.4) \qquad y_{ij} = \max(0, x_{ij}^T \alpha_0 + u_{ij}),$$

versus the local alternative model

$$(2.5) \qquad y_{ij} = \max(0, x_{ij}^T \alpha_0 + n^{-1/2} z_{ij}^T \beta_0 + u_{ij}).$$

To derive the quantile estimate of $\alpha_0$ in (2.4), we follow Powell [25] and consider the minimization of the objective function

$$Q_n(\alpha) = \sum_{ij} \rho_\tau \{y_{ij} - \max(0, x_{ij}^T \alpha)\}, \tag{2.6}$$

where $\rho_\tau(u) = u \cdot \{\tau - I(u < 0)\}$ is the quantile loss function. Under mild conditions, it is established in Theorem 2.1 below that the quantile estimator $\hat{\alpha}_0$ is strongly consistent and asymptotically normal in both models (2.4) and (2.5) even though the objective function (2.6) treats all observations as if they were independent.

In the absence of censored observations, minimization of (2.6) can be performed efficiently by linear programming techniques. In fact, the solution to (2.6) in such cases only requires software for quantile regression in linear models. However, with censored observations, the objective function in (2.6) is neither differentiable nor convex, and this presents a computational challenge. A number of optimization approaches have been proposed in the literature; see, for example, [1, 6, 20]. In this paper we employ the BRECNS algorithm of [6] as implemented in the R package *quantreg*.

2.2. *Large sample properties of $\hat{\alpha}_0$.* Throughout the paper we suppose a typical longitudinal data set where $n_i$ (the number of repeated measurements for each subject) is bounded, but $N$ (the number of subjects) grows. Note that all results are stated for a given $\tau$, although dependence on $\tau$ is not explicit in the various expressions. To establish all the large sample properties in this paper, we require the following conditions:

A1. The parameter vector $\alpha_0$ is an interior point of a compact parameter space $\mathbb{A} \in \mathbb{R}^p$.
A2. Let $\|\tilde{x}_{ij}\|$ denote the Euclidean norm of $\tilde{x}_{ij}$, where $\tilde{x}_{ij} = (x_{ij}^T, z_{ij}^T)^T$, then $\max_{ij} \|\tilde{x}_{ij}\| = O(n^{1/4})$ and $n^{-1} \sum_{ij} \|\tilde{x}_{ij}\|^3 = O(1)$ as $n \to \infty$.
A3. There exists $\varepsilon_0 > 0$ such that as $n \to \infty$, $\liminf n^{-1} \sum_{ij} I(|x_{ij}^T \alpha| \geq \varepsilon_0) > 0$ for any $\|\alpha\| \neq 0$ and $D_{1n}(\alpha_0) = n^{-1} \sum_{ij} I(x_{ij}^T \alpha_0 \geq \varepsilon_0) x_{ij} x_{ij}^T \to D_1$, where $D_1$ is a positive definite matrix.
A4. The $u_{ij}$ have a common marginal distribution function $F$ and a Lebesque density $f$, which is Lipschitz in a neighborhood of 0. Also, there exist some positive values $\varrho_1$ and $\varrho_2$ such that $f(u) < \varrho_2$ for all $u$, and $f(u) > \varrho_1$ for $|u| < \varrho_1$.
A5. For any $d \geq 0$, there exists a positive constant $C$ such that $n^{-1} \sum_{ij} I(|x_{ij}^T \alpha_0| \leq \|x_{ij}\| d) \leq Cd$.
A6. Let $D_{2n}(\alpha_0) = n^{-1} \sum_{ij} I(x_{ij}^T \alpha_0 > 0) x_{ij} z_{ij}^T \to D_2$, as $n \to \infty$, where $D_2$ is a $p \times q$ matrix.
A7. The joint distribution function of $u_{ij_1}$ and $u_{ij_2}$ for any $i$ and $j_1 \neq j_2$, denoted as $F_{1,2}$, is Lipschitz in a neighborhood of $(0,0)$.

A8. Let $D_{3n}(\alpha_0) = n^{-1} \sum_{ij} I(x_{ij}^T \alpha_0 > 0) z_{ij}^* z_{ij}^T \to D_3$, as $n \to \infty$.

The following theorem states the large sample properties of $\hat{\alpha}_0$ under models (2.4) and (2.5):

THEOREM 2.1. *For the longitudinal censored regression models (2.4) and (2.5):*

(i) *If conditions* A1–A4 *hold, then the censored quantile estimator $\hat{\alpha}_0$ converges to $\alpha_0$ almost surely.*

(ii) *If conditions* A1–A5 *hold, then under model (2.4) the censored quantile estimator $\hat{\alpha}_0$ is asymptotically normal,*

$$\{\Gamma_n(\delta)\}^{-1/2} \sqrt{n}(\hat{\alpha}_0 - \alpha_0) \xrightarrow{D} N(0, I).$$

(iii) *If conditions* A1–A6 *hold, then under model (2.5) the censored quantile estimator $\hat{\alpha}_0$ is asymptotically normal,*

$$\{\Gamma_n(\delta)\}^{-1/2} \sqrt{n}(\hat{\alpha}_0 - \alpha_0 - D_1^{-1} D_2 \beta_0) \xrightarrow{D} N(0, I),$$

*where*

$$\Gamma_n(\delta) = n^{-1} \{f(0)\}^{-2} D_1^{-1}$$
$$\times \left\{ \sum_{ij} I(x_{ij}^T \alpha_0 > 0) x_{ij} x_{ij}^T \tau(1-\tau) \right.$$
$$\left. + \sum_i \sum_{j \neq j'} I(x_{ij}^T \alpha_0 > 0, x_{ij'}^T \alpha_0 > 0) x_{ij} x_{ij}^T (-\tau^2 + \delta) \right\} D_1^{-1}$$

*and $\delta = P(u_{i1} < 0, u_{i2} < 0)$ measures the intra-subject dependence.*

REMARK 1. The common density assumption of $u_{ij}$ in A4 is made for convenience, but not necessary for the strong consistency nor the asymptotic normality of $\hat{\alpha}_0$. In order for Theorem 2.1 to hold, it suffices that the $\tau$th quantile of $u_{ij}$ is 0 for all $i$ and $j$ with density functions $f_{ij}$, which are continuously differentiable in a neighborhood of zero and uniformly bounded away from zero and infinity. When $y_{ij}$ is left censored at some known values $c_{ij}$, the asymptotic results developed in this paper hold, but $x_{ij}^T \alpha_0$ in conditions A3, A5, A6 and A8 must be replaced with $x_{ij}^T \alpha_0 - c_{ij}$.

2.3. *Quantile rank score test.* To test the hypotheses in (2.3), one can explore the asymptotic normality of censored quantile estimators of the parameters $(\alpha_0, \beta)$ in model (2.2). However, following the proof of Theorem 2.1 part (ii), one can see that the asymptotic variance–covariance matrix of these estimators is a function of the unspecified density of error terms.

This hampers the use of a Wald-type test. Moreover, it has been shown that, in a quantile regression set up, a Wald-type test is generally unstable at small sample sizes (e.g., [3, 18]). The use of likelihood ratio-based tests is even more daunting for our setup because the limiting distribution is a complicated function of the unknown error density. To avoid these problems and the need for estimating a density, which is in our testing problem an infinite dimensional nuisance parameter, we turn to the quantile rank score test proposed in [9] for independent and uncensored data.

To present our test, we rewrite model (2.2) in matrix form

$$(2.7) \qquad Y = \max(\underline{0}_n, X\alpha_0 + Z\beta + U),$$

where $Y$ and $U$ are $n$-dimensional vectors, $\underline{0}_n$ is an $n \times 1$ vector consisting of zeros and $X$ and $Z$ are $n \times p$ and $n \times q$ matrices, respectively. Let $X^* = \text{diag}\{I(X\alpha_0 > 0)\}X$, $H = X^*(X^{*T}X^*)^{-1}X^{*T}$ and $Z^* = (z_{ij}^*)_{n \times q} = (I - H)Z$. Note that $Z^*$, which is a linear combination of the design matrix $Z$, is orthogonal to the space spanned by those $x_{ij}$'s that satisfy $x_{ij}^T\alpha_0 > 0$.

Our proposed quantile rank score test is based on

$$(2.8) \qquad S_n = n^{-1/2} \sum_{ij} \{I(x_{ij}^T\hat{\alpha}_0 > 0)z_{ij}^*\varphi_\tau(\hat{u}_{ij})\},$$

where $\hat{u}_{ij} = y_{ij} - \max(0, x_{ij}^T\hat{\alpha}_0)$, $\hat{\alpha}_0$ is the censored quantile estimator of $\alpha_0$ in model (2.4), and $\varphi_\tau(u) = \tau - I(u < 0)$ is the quantile score function. It is worth pointing out that $\varphi_\tau(u)$ is the piecewise gradient of the quantile loss function $\rho_\tau(u)$, and that $S_n$ only includes scores from those observations for which the corresponding $x_{ij}^T\alpha_0$ are estimated to be uncensored.

Let

$$(2.9) \qquad V_n(\delta; \hat{\alpha}_0) = n^{-1}\Bigg\{\sum_{ij} I(x_{ij}^T\hat{\alpha}_0 > 0)z_{ij}^*z_{ij}^{*T}\tau(1-\tau) \\ + \sum_{i,j \neq j'} I(x_{ij}^T\hat{\alpha}_0 > 0, x_{ij'}^T\hat{\alpha}_0 > 0)z_{ij}^*z_{ij'}^{*T}(-\tau^2 + \delta)\Bigg\},$$

where $\delta$ is defined in Theorem 2.1(iii).

We define the Quantile Rank Score ($QRS$) test statistic as

$$(2.10) \qquad T_n = S_n^T\{V_n(\hat{\delta}; \hat{\alpha}_0)\}^{-1}S_n,$$

where $\hat{\delta} = L^{-1}\sum_{i,j \neq j'} I(x_{ij}^T\hat{\alpha}_0 > 0, x_{ij'}^T\hat{\alpha}_0 > 0)I(\hat{u}_{ij} < 0, \hat{u}_{ij'} < 0)$ and $L$ denotes the total number of pairs of repeated measurements that are predicted to be uncensored. Note that when all the observations are uncensored and independent (i.e., $\delta = \tau^2$ and $n_i = 1$), $T_n$ reduces to the $QRS$ test-statistic proposed in [9].

THEOREM 2.2. *Assume that conditions* A1–A8 *hold, then as* $n \to \infty$, *we have:*

(i) *under* $H_0$, *the statistic* $T_n$ *is asymptotically* $\chi^2$ *with* $q$ *degrees of freedom;*

(ii) *under* $H_n$, $T_n$ *is asymptotically noncentral* $\chi^2$ *with* $q$ *degrees of freedom and with noncentrality parameter* $\beta_0^T D_3 [V_n(\delta; \alpha_0)]^{-1} D_3 \beta_0 f^2(0)$.

REMARK 2. The joint probability $\delta$ in (2.9) captures the sign correlation between errors from the same subject. When $\delta \in (\tau^2, \tau]$ these errors are positively correlated, when $\delta \in [0, \tau^2)$ they are negatively correlated and when $\delta = \tau^2$ the errors are independent. Ignoring the intra-subject dependence leads to a test $T_n^*$, say. Depending on $\delta$ and/or $Z$, this test statistic is either invalid or lacks power. For illustration, consider the case where $q = 1$ and $\delta \in (\tau^2, \tau]$. Then we have

$$(2.11) \quad \frac{V_n(\delta; \alpha_0)}{V_n(\tau^2; \alpha_0)} = 1 + \frac{\sum_{i=1} \sum_{j \neq j'} I(x_{ij}^T \alpha_0 > 0, x_{ij'}^T \alpha_0 > 0) z_{ij}^* z_{ij'}^* (\delta - \tau^2)}{\sum_{ij} I(x_{ij}^T \alpha_0 > 0) z_{ij}^{*2} (\tau - \tau^2)}.$$

When testing the between-subject factor effect, for example in model (3.1) of the simulation study, the signs of $z_{ij}^* z_{ij'}^*$ are positive for the same $i$th subject, $\frac{V_n(\delta; \alpha_0)}{V_n(\tau^2; \alpha_0)} > 1$ and $T_n^*$ leads to inflated Type I errors. For a given significance level $\theta$, the power of $T_n^*$ under $H_n$ is

$$1 - \Phi \left\{ Z_{\theta/2} \sqrt{\frac{V_n(\tau^2; \alpha_0)}{V_n(\delta; \alpha_0)}} - \frac{\mu_n}{\sqrt{V_n(\delta; \alpha_0)}} \right\}$$
$$+ \Phi \left\{ -Z_{\theta/2} \sqrt{\frac{V_n(\tau^2; \alpha_0)}{V_n(\delta; \alpha_0)}} - \frac{\mu_n}{\sqrt{V_n(\delta; \alpha_0)}} \right\},$$

where $\mu_n = n^{-1} f(0) \sum_{ij} I(x_{ij}^T \alpha_0 > 0) z_{ij}^* z_{ij} \beta_0$, $\Phi$ denotes the CDF of the standard normal distribution and $Z_\theta$ is the upper $\theta$th quantile of $\Phi$. Therefore, when testing the within-subject factor effect, for example, we can have the signs of all $z_{ij}^* z_{ij'}^*$ be negative for the same subject, thus $\frac{V_n(\delta; \alpha_0)}{V_n(\tau^2; \alpha_0)} < 1$ and $T_n^*$ has diminished power.

2.4. *Construction of confidence intervals.*

2.4.1. *Confidence intervals via inversion of rank score tests.* The developed rank score test can be extended to test $H_0 : \beta = \beta_0$ by simply rewriting model (2.7). We denote $\tilde{y}_{ij} = y_{ij} - z_{ij}^T \beta_0$. The fact that $y_{ij}$ is censored at 0 implies that $\tilde{y}_{ij}$ is censored from the left at $-z_{ij}^T \beta_0$. It is clear that under $H_0$, the $\tau$th quantile estimate of $\alpha$ can be obtained by minimizing $\sum_{ij} \rho_\tau \{\tilde{y}_{ij} - \max(-z_{ij}^T \beta_0, x_{ij}^T \alpha)\}$. The quantile rank score test can be constructed following the same procedure as in Section 2.3 by replacing $y$ with $\tilde{y}$.

For a quantile coefficient $\beta \in \mathbb{R}^1$, the confidence interval can be constructed by inverting the rank score test. Using the fact that the test statistic $T_n$ is convex in $\beta$, we can obtain a $100(1-\theta)\%$ confidence interval consisting of the $\beta_0$'s at which the test on $H_0: \beta = \beta_0$ will not be rejected. The readers are referred to [18, 19] or [3] for details of confidence interval construction for uncensored and independent data.

2.4.2. *Blockwise modified bootstrap method.* The more computationally demanding resampling method offers an alternative approach for statistical inference. Here we introduce a modified bootstrap approach through blockwise pairs resampling, denoted by *Boot*. For easy presentation, we denote $\gamma = (\alpha_0^T, \beta^T)^T \in \mathbb{R}^{p+q}$.

In applications of quantile regression, the pairs bootstrap is often chosen over the residual bootstrap because it is insensitive to model misspecification and heteroscedasticity. The idea of the pairs bootstrap is to draw pairs, in our case, $(y_{ij}, \tilde{x}_{ij})$, at random from the original observations with replacement. Note that in model (2.1), the observations are dependent within each subject. To retain this dependence structure, we treat the observations in each subject as a block and resample the block pairs $\{(y_{ij}, \tilde{x}_{ij}), j = 1, \ldots, n_i\}$.

As we have seen, computation of Powell's estimator is complicated by the nonconvexity of the objective function (2.6). Therefore, direct implementation of the bootstrap approach could be prohibitively expensive in terms of computation. To reduce the computational cost, we employ a modified bootstrap method proposed in [2] for median regression with independent data. From now on, we define the bootstrap sample of $(y_{ij}, \tilde{x}_{ij})$ by $(y_{ij}^\#, x_{ij}^\#)$. It is known that the solution of $\min_{a \in \mathbb{R}^p} \sum_{ij} \rho_\tau(y_{ij} - \tilde{x}_{ij}^T \gamma) I(\tilde{x}_{ij}^T \gamma_0 > 0)$ is asymptotically equivalent to Powell's estimator. Making use of the $\hat{\gamma}$ that result from fitting the model with the observed data, the modified bootstrap estimator $\hat{\gamma}^\#$ can be obtained by minimizing

$$(2.12) \qquad \sum_{ij} \rho_\tau(y_{ij}^\# - x_{ij}^{\#T} \gamma) I(x_{ij}^{\#T} \hat{\gamma} > 0).$$

Note that (2.12) is a convex function, and thus $\hat{\gamma}^\#$ can be calculated in the same way as in uncensored quantile regression. A $100(1-\theta)\%$ confidence interval for $\gamma$ can be obtained with the lower and upper bound calculated as the $(\theta/2)$th and $(1-\theta/2)$th quantiles of those $\hat{\gamma}^\#$'s.

**3. Simulation study.** To assess the performance of the inference procedures described in Section 2, we conduct a simulation study. We explore the effects of different proportions of censoring and various degrees of intrasubject dependency on estimation, testing and confidence intervals.

3.1. *Model descriptions.* In the simulation, the latent response variable $y^*$ is generated from the following model:

$$(3.1) \quad y_{ij}^* = 1 + x_{ij}\alpha + z_{ij}\beta + \sigma_{ij}\{u_{ij} - F_u^{-1}(\tau)\},$$
$$i = 1, \ldots, N, \ j = 1, \ldots, 10,$$

where $u_{ij} = a_i + e_{ij}$ is the random error and $F_u^{-1}(\tau)$ is the $\tau$th quantile of $u$, $a_i$ is the random subject effect, $x_{ij}$ and $e_{ij}$ are i.i.d. from the standard normal distribution, $z_{ij} = 0$ for the first $N/2$ subjects and $z_{ij} = 1$ for the rest. Four different cases are considered:

Case 1. A fixed effect model ($a_i = 0$) with homoscedastic term $\sigma_{ij} = 1$.

Case 2. A random effect model with $a_i$ that are i.i.d. from the standard normal and $\sigma_{ij} = 1$. This yields a homoscedastic model with an intra-subject correlation coefficient of 0.5.

Case 3. A random effect model with $a_i$ that are i.i.d. from $N(0,9)$ and $\sigma_{ij} = 1$. This yields a homoscedastic model with an intra-subject correlation coefficient of 0.9.

Case 4. A heteroscedastic model with $a_i \sim N(0,1)$ and $\sigma_{ij} = 1 + |x_{ij}|$.

For all cases we consider both 20% and 40% censoring. The observed response variable $y_{ij}$ is generated from the maximum of 0 and $y_{ij}^*$, subtracting the 20th or 40th percentile of $\{y_{ij}^*\}$, respectively. Our analysis focuses on the effect of $z_{ij}$ at three quartiles, as per the main objective of Section 2. To evaluate Type I error and power, the nominal significance level is set to 5%, $\alpha$ is fixed at 10, $\beta$ is varied from 0 to 1 and the simulation was repeated 500 times in all cases.

3.2. *Evaluation of the proposed estimator for $\alpha$.* We first compare the finite sample efficiency of the omniscient estimator, our estimator $\hat{\alpha}$ and two naïve estimators. The omniscient estimator is obtained by fitting the quantile regression model with the latent response variable and thus serves as a gold standard. The first naïve estimator, *Naïve*$_1$, is obtained using all observations as if none were censored. The second, *Naïve*$_2$, is computed using only the uncensored observations.

Table 1 summarizes the bias and mean squared error (MSE) of these four estimators in Cases 1–4 for $N = 50$ and $\beta = 0$. Compared to the omniscient estimator, $\hat{\alpha}$ performs universally well, even when the data are highly correlated (Case 3). As expected, the two naïve estimators have larger biases and mean squared errors than $\hat{\alpha}$, even more so for the higher proportion of censoring.

TABLE 1
*Comparison of the omniscient estimator, our estimator $\hat{\alpha}$ and two naïve estimators for $\alpha$ at $N = 50$. The Naïve$_1$ is obtained by uncensored quantile regression and Naïve$_2$ is obtained using only the uncensored observations. The CP stands for the censoring proportion and MSE stands for mean squared error*

|  |  | Omniscient |  | $\hat{\alpha}$ |  | Naïve$_1$ |  | Naïve$_2$ |  |
| --- | --- | --- | --- | --- | --- | --- | --- | --- | --- |
|  | CP | MSE | Bias | MSE | Bias | MSE | Bias | MSE | Bias |
| | | | | $\tau = 0.25$ | | | | | |
| Case 1 | 0.2 | 0.004 | 0 | 0.008 | 0.008 | 0.360 | $-0.591$ | 0.013 | $-0.072$ |
|  | 0.4 | 0.004 | 0 | 0.015 | 0.014 | 6.791 | $-2.576$ | 0.033 | $-0.135$ |
| Case 2 | 0.2 | 0.007 | $-0.007$ | 0.017 | 0.001 | 0.719 | $-0.834$ | 0.043 | $-0.159$ |
|  | 0.4 | 0.007 | $-0.007$ | 0.032 | 0.013 | 9.408 | $-3.039$ | 0.119 | $-0.295$ |
| Case 3 | 0.2 | 0.037 | $-0.013$ | 0.093 | 0.005 | 2.953 | $-1.685$ | 0.699 | $-0.762$ |
|  | 0.4 | 0.037 | $-0.013$ | 0.180 | 0.019 | 20.953 | $-4.548$ | 2.266 | $-1.427$ |
| Case 4 | 0.2 | 0.043 | $-0.014$ | 0.098 | 0.017 | 2.621 | $-1.594$ | 0.392 | $-0.550$ |
|  | 0.4 | 0.043 | $-0.014$ | 0.183 | 0.051 | 16.572 | $-4.045$ | 0.584 | $-0.647$ |
| | | | | $\tau = 0.5$ | | | | | |
| Case 1 | 0.2 | 0.003 | 0.001 | 0.007 | 0.006 | 0.717 | $-0.838$ | 0.012 | $-0.067$ |
|  | 0.4 | 0.003 | 0.001 | 0.014 | 0.004 | 10.878 | $-3.277$ | 0.033 | $-0.136$ |
| Case 2 | 0.2 | 0.006 | $-0.001$ | 0.013 | 0.002 | 1.142 | $-1.058$ | 0.031 | $-0.138$ |
|  | 0.4 | 0.006 | $-0.001$ | 0.023 | 0.007 | 11.686 | $-3.401$ | 0.094 | $-0.267$ |
| Case 3 | 0.2 | 0.032 | $-0.009$ | 0.063 | $-0.006$ | 2.782 | $-1.646$ | 0.448 | $-0.610$ |
|  | 0.4 | 0.032 | $-0.009$ | 0.115 | 0.017 | 15.583 | $-3.931$ | 1.434 | $-1.133$ |
| Case 4 | 0.2 | 0.021 | 0.002 | 0.040 | 0.009 | 2.218 | $-1.476$ | 0.129 | $-0.304$ |
|  | 0.4 | 0.035 | $-0.002$ | 0.128 | 0.035 | 17.266 | $-4.141$ | 0.436 | $-0.568$ |
| | | | | $\tau = 0.75$ | | | | | |
| Case 1 | 0.2 | 0.004 | 0.002 | 0.008 | 0.006 | 2.408 | $-1.537$ | 0.012 | $-0.062$ |
|  | 0.4 | 0.004 | 0.002 | 0.015 | 0.002 | 15.134 | $-3.879$ | 0.030 | $-0.123$ |
| Case 2 | 0.2 | 0.007 | $-0$ | 0.013 | 0.005 | 2.724 | $-1.635$ | 0.027 | $-0.118$ |
|  | 0.4 | 0.007 | $-0$ | 0.023 | 0.008 | 15.157 | $-3.883$ | 0.080 | $-0.237$ |
| Case 3 | 0.2 | 0.038 | $-0.008$ | 0.063 | 0.005 | 3.384 | $-1.818$ | 0.340 | $-0.514$ |
|  | 0.4 | 0.038 | $-0.008$ | 0.101 | $-0.009$ | 14.717 | $-3.824$ | 0.998 | $-0.927$ |
| Case 4 | 0.2 | 0.040 | $-0.001$ | 0.058 | 0.010 | 4.420 | $-2.085$ | 0.222 | $-0.403$ |
|  | 0.4 | 0.040 | $-0.001$ | 0.116 | 0.019 | 20.051 | $-4.467$ | 0.372 | $-0.509$ |

3.3. *Performance of the proposed quantile rank score test.* We evaluate the performance of our proposed Quantile Rank Score test ($QRS$) by comparing it to five other test statistics: a rank score test for censored data that assumes independence (*Indep*); a naïve rank score test (*Naïve$_1$*) that assumes the observations are uncensored; another naïve rank score test (*Naïve$_2$*) that uses only the uncensored observations; a bootstrap base test, *Boot*, with 500 resamplings; and the omniscient rank score test, *Omni*, based on the latent response variable.

TABLE 2
*The Type* I *errors in Cases 1–4 at* $N = 10$ *and* $N = 50$

|  | Case 1 | | | Case 2 | | | Case 3 | | | Case 4 | | |
|---|---|---|---|---|---|---|---|---|---|---|---|---|
| $\tau$ | 0.25 | 0.5 | 0.75 | 0.25 | 0.5 | 0.75 | 0.25 | 0.5 | 0.75 | 0.25 | 0.5 | 0.75 |
| | | | | | $CP = 0.2, N = 10$ | | | | | | | |
| *Indep* | 0.04 | 0.06 | 0.06 | 0.30 | 0.30 | 0.26 | 0.42 | 0.43 | 0.42 | 0.30 | 0.22 | 0.27 |
| *Naïve*$_1$ | 0.05 | 0.05 | 0.04 | 0.07 | 0.07 | 0.05 | 0.07 | 0.07 | 0.05 | 0.07 | 0.06 | 0.04 |
| *Naïve*$_2$ | 0.04 | 0.08 | 0.09 | 0.07 | 0.07 | 0.09 | 0.09 | 0.09 | 0.13 | 0.07 | 0.08 | 0.12 |
| *Omni* | 0.05 | 0.05 | 0.06 | 0.06 | 0.06 | 0.05 | 0.07 | 0.06 | 0.06 | 0.05 | 0.05 | 0.06 |
| *QRS* | 0.05 | 0.05 | 0.06 | 0.06 | 0.05 | 0.05 | 0.07 | 0.06 | 0.04 | 0.06 | 0.05 | 0.06 |
| *Boot* | 0.05 | 0.07 | 0.07 | 0.09 | 0.08 | 0.09 | 0.09 | 0.07 | 0.10 | 0.10 | 0.08 | 0.04 |
| | | | | | $CP = 0.4, N = 10$ | | | | | | | |
| *Indep* | 0.05 | 0.07 | 0.05 | 0.26 | 0.24 | 0.23 | 0.37 | 0.39 | 0.38 | 0.26 | 0.19 | 0.22 |
| *Naïve*$_1$ | 0.04 | 0.06 | 0.04 | 0.05 | 0.04 | 0.04 | 0.06 | 0.06 | 0.05 | 0.06 | 0.05 | 0.05 |
| *Naïve*$_2$ | 0.05 | 0.09 | 0.10 | 0.10 | 0.09 | 0.12 | 0.12 | 0.14 | 0.23 | 0.11 | 0.12 | 0.13 |
| *Omni* | 0.05 | 0.05 | 0.06 | 0.06 | 0.06 | 0.05 | 0.07 | 0.06 | 0.06 | 0.05 | 0.05 | 0.06 |
| *QRS* | 0.04 | 0.08 | 0.05 | 0.06 | 0.07 | 0.06 | 0.07 | 0.07 | 0.07 | 0.06 | 0.07 | 0.05 |
| *Boot* | 0.05 | 0.08 | 0.06 | 0.08 | 0.07 | 0.10 | 0.10 | 0.08 | 0.13 | 0.09 | 0.10 | 0.10 |
| | | | | | $CP = 0.2, N = 50$ | | | | | | | |
| *Indep* | 0.05 | 0.06 | 0.04 | 0.29 | 0.27 | 0.29 | 0.44 | 0.40 | 0.43 | 0.28 | 0.20 | 0.28 |
| *Naïve*$_1$ | 0.05 | 0.04 | 0.04 | 0.06 | 0.04 | 0.06 | 0.05 | 0.05 | 0.03 | 0.05 | 0.05 | 0.05 |
| *Naïve*$_2$ | 0.06 | 0.06 | 0.06 | 0.05 | 0.05 | 0.06 | 0.05 | 0.06 | 0.08 | 0.05 | 0.06 | 0.06 |
| *Omni* | 0.05 | 0.05 | 0.04 | 0.04 | 0.05 | 0.04 | 0.05 | 0.04 | 0.05 | 0.04 | 0.06 | 0.04 |
| *QRS* | 0.05 | 0.06 | 0.05 | 0.04 | 0.04 | 0.04 | 0.06 | 0.04 | 0.04 | 0.04 | 0.06 | 0.04 |
| *Boot* | 0.04 | 0.06 | 0.04 | 0.05 | 0.05 | 0.06 | 0.06 | 0.04 | 0.06 | 0.05 | 0.06 | 0.05 |
| | | | | | $CP = 0.4, N = 50$ | | | | | | | |
| *Indep* | 0.06 | 0.05 | 0.05 | 0.22 | 0.24 | 0.22 | 0.37 | 0.33 | 0.35 | 0.23 | 0.17 | 0.23 |
| *Naïve*$_1$ | 0.05 | 0.04 | 0.06 | 0.04 | 0.04 | 0.08 | 0.05 | 0.04 | 0.06 | 0.06 | 0.06 | 0.05 |
| *Naïve*$_2$ | 0.05 | 0.06 | 0.06 | 0.06 | 0.05 | 0.08 | 0.07 | 0.08 | 0.10 | 0.05 | 0.06 | 0.09 |
| *Omni* | 0.05 | 0.05 | 0.04 | 0.04 | 0.05 | 0.04 | 0.05 | 0.04 | 0.05 | 0.04 | 0.06 | 0.04 |
| *QRS* | 0.05 | 0.05 | 0.04 | 0.04 | 0.04 | 0.04 | 0.06 | 0.04 | 0.04 | 0.04 | 0.04 | 0.04 |
| *Boot* | 0.04 | 0.06 | 0.05 | 0.04 | 0.04 | 0.05 | 0.06 | 0.06 | 0.07 | 0.03 | 0.05 | 0.05 |

Table 2 summarizes the Type I error rates of all six test statistics in Cases 1–4 for $N = 10$ or 50. The Type I error of *Boot* is estimated by the proportion of cases where 0 is not contained in the 95% confidence interval for $\beta$. It is obvious that the rank score test *Indep* without the $\delta$ adjustment loses complete control of Type I error in Cases 2, 3 and 4, where the data are correlated due to the random effects $a_i$. Moreover, the size of the deterioration increases with the degree of intra-subject dependency. As to the naïve tests, we find that, in general, *Naïve*$_2$ has inflated Type I errors at $N = 10$, and *Naïve*$_1$ lacks power (see Figure 1). The modified bootstrap method, *Boot*, preserves the nominal significance level well at $N = 50$, but gives consistently inflated Type I errors in the smaller samples $N = 10$. Generally speaking, the *QRS* and omniscient methods preserve the nominal significance level reasonably well in all cases.

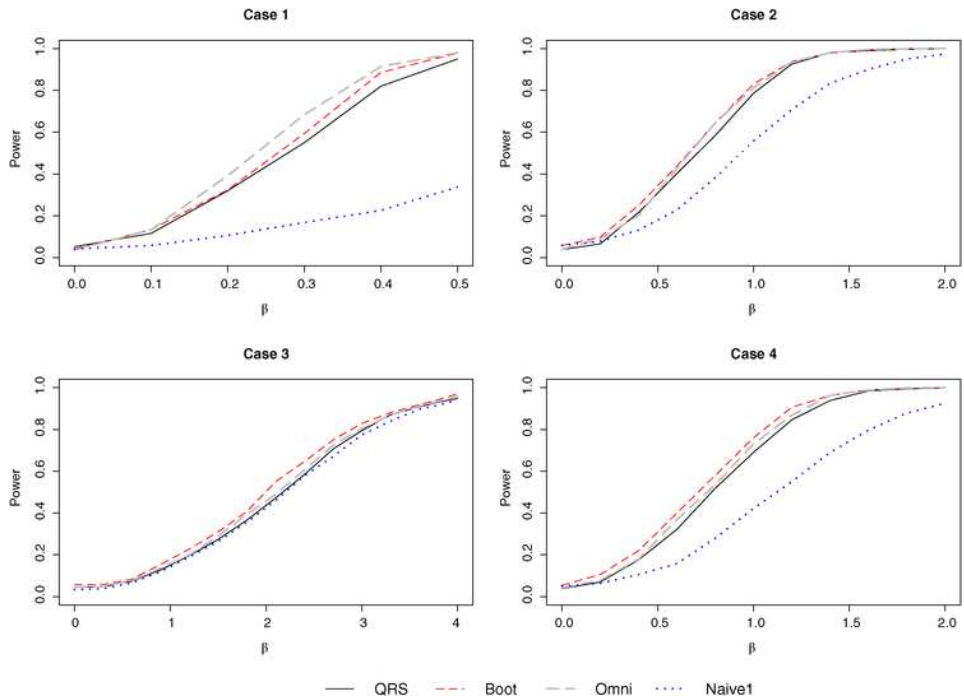

FIG. 1. *Power curves of QRS, Boot, Omni and Naïve$_1$ in Cases 1–4 at $\tau = 0.75$ with $N = 50$ and 20% censoring.*

Figure 1 plots the power curves of *QRS*, *Boot*, *Omni* and *Naïve$_1$* at $\tau = 0.75$ with $N = 50$ and 20% censoring. The *Naïve$_1$* loses a great deal of power by ignoring the censoring. The *QRS* and *Boot* both perform as well as the omniscient method in all cases. Using the modified bootstrap approach reduces the computational time as compared to the direct implementation of bootstrap, but it is still much more computationally intensive than *QRS*. For example, using R (version 2.3.1) in a 3.4 GHz Dell computer with 3.0 GB of RAM to simulate Case 2 at $\tau = 0.75$ with $N = 10$ and 20% censoring, the *QRS* took 23 seconds for 500 runs of simulation, compared to 3,013 seconds by *Boot* with 500 resamples. Furthermore, *QRS* is robust to the heteroscedasticity considered in Case 4, even though it is developed for models with homoscedastic errors.

We study the performance of *QRS* under the local alternative $H_n : \beta = n^{-1/2}\beta_0$ for $n = 10N$ varying from 200 to 5000. We let $\beta_0 = 10$ in Case 1, $\beta_0 = 20$ in Cases 2 and 4 and $\beta_0 = 50$ in Case 3. Figure 2 shows that the local power of *QRS* remains stable as $n$ increases. This observation is consistent with the asymptotic results in Section 2.3. Note that *QRS* exhibits different powers in four cases due to the different variation used to generate the subject effects.

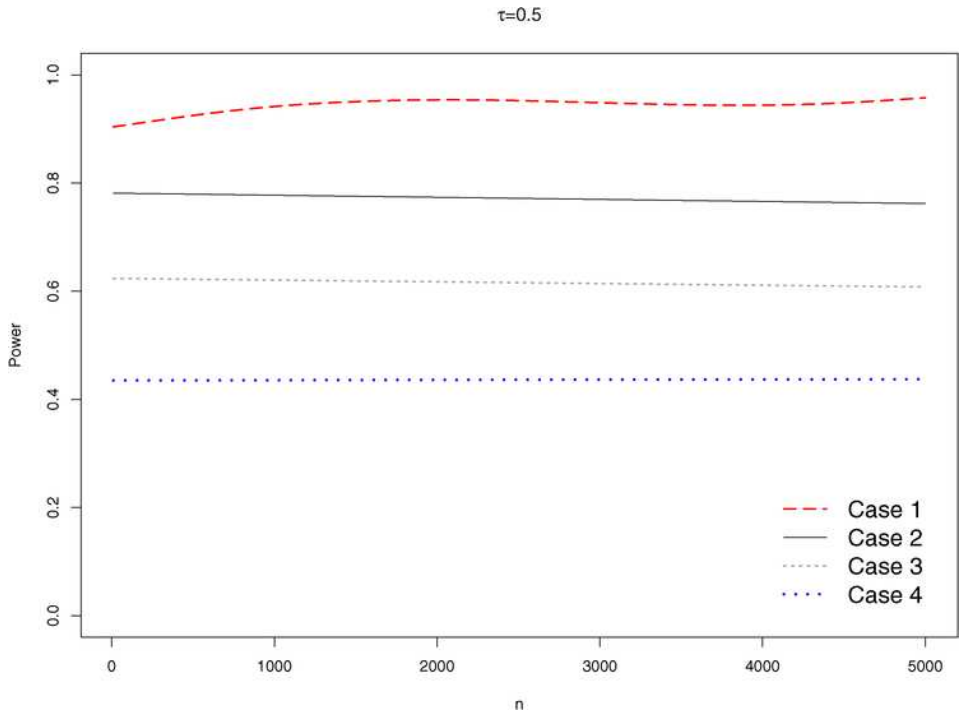

FIG. 2.   *Power curves of QRS under local alternative $H_n: \beta = 50n^{-1/2}$ in Cases 1–4 at $\tau = 0.5$ with 20% censoring. All the curves are generated by fitting smoothing splines over the estimated powers against $n$.*

3.3.1. *Assessment of confidence intervals for $\beta$.* For each simulated data set under the hypothesis $\beta = 1$, we obtain 95% confidence intervals for $\beta$ using the bootstrap method, *Boot*, and by inverting the *QRS* test following the procedure described in Section 2.4. For comparison, we also obtain confidence intervals by inverting the *Omni* test and the two naïve tests, *Naïve*$_1$ and *Naïve*$_2$. The estimated mean lengths (EML) of these five confidence interval procedures and the empirical coverage probabilities (ECP) across the 500 intervals in Case 2 at $N = 50$ are summarized in Table 3.

The two naïve methods give poor confidence intervals at all quantile levels. The empirical coverage probabilities of *Omni*, *QRS* and *Boot* are, in general, close to the nominal level. Moreover, the mean lengths of the *QRS* confidence intervals are comparable to those of the *Omni*.

**4. Application to an HIV-RNA level study.** In this section we apply the methodology proposed in Section 2 to analyze the HIV-RNA data in [10]. This clinical trial followed a total of 481 HIV-infected individuals with baseline HIV-RNA levels in their plasma (viral load) greater than 1000 copies/ml. For each individual, viral load was measured at time zero and then approximately 2, 4, 8, 16 and 24 weeks later. Due to the detection limit of the assay used to measure viral load, 22% of measurements were censored from below at 200 copies/ml. We refer the readers to [10] for a more thorough discussion of this AIDS-related clinical trail. We seek to compare the viral load response (VL) to a double protease inhibitor (DPI) regimen, herein referred to as the treatment, with that of a single protease inhibitor (SPI) regimen, herein referred to as the control.

Our preliminary investigation shows that the $\log_{10}$VL from both regimens drops sharply during the first two weeks, and then remains rather stable. To capture this evolution pattern, we consider a model that includes an intercept, a slope for the first two weeks, and a different slope for the remainder of the study. The working model is

$$(4.1) \quad \begin{aligned} y_{ij} = \max\{\log_{10}(200), &\beta_0 + \beta_1 \min(t_{ij}, 2) + \beta_2(t_{ij} - 2)I(t_{ij} > 2) \\ &+ \gamma_1 \min(t_{ij}, 2)z_i + \gamma_2(t_{ij} - 2)I(t_{ij} > 2)z_i + u_{ij}\}, \end{aligned}$$

where $y_{ij}$ is the observed $\log_{10}$VL of the $i$th subject at time $t_{ij}$, $z_i$ is the treatment indicator taking 1 for the treatment group and 0 for the control.

We fit model (4.1) at different quantiles with $\tau$ varying from 0.25 to 0.9 to obtain a profile of the regimens' effects. Sun and Wu [28] analyzed the same data using a partial linear model. However, they ignored the left censored observations and used only those responses within the detectable range. We shall see that doing so leads to biased results and distorted conclusions about the treatment's efficacy.

Figure 3 shows the estimated quartiles of two regimens at each time point from our method (curves with solid points), and those from Naïve$_2$ that use only the uncensored observations. As expected, the naïve method overestimates VL, especially in the lower quantiles. At $\tau = 0.25$, our method shows that the VL drops rapidly during the first two weeks, and it continues dropping afterward for both regimens, while *Naïve$_2$* suggests an increasing trend after week 2.

TABLE 3
*The empirical coverage probabilities (ECP) and estimated mean lengths (EML) for confidence interval procedures in Case 2 at $N = 50$. The nominal level is 0.95*

| Method | $\tau = 0.25$ | | $\tau = 0.50$ | | $\tau = 0.75$ | |
|---|---|---|---|---|---|---|
| | ECP | Length | ECP | Length | ECP | Length |
| *Omni* | 0.95 | 0.79 | 0.94 | 0.74 | 0.96 | 0.78 |
| *QRS* | 0.96 | 1.00 | 0.93 | 0.79 | 0.95 | 0.80 |
| *Boot* | 0.95 | 1.26 | 0.95 | 1.26 | 0.95 | 1.26 |
| *Naïve$_1$* | 0.34 | 0.56 | 0.77 | 0.63 | 0.92 | 0.72 |
| *Naïve$_2$* | 0.28 | 0.54 | 0.59 | 0.55 | 0.71 | 0.56 |

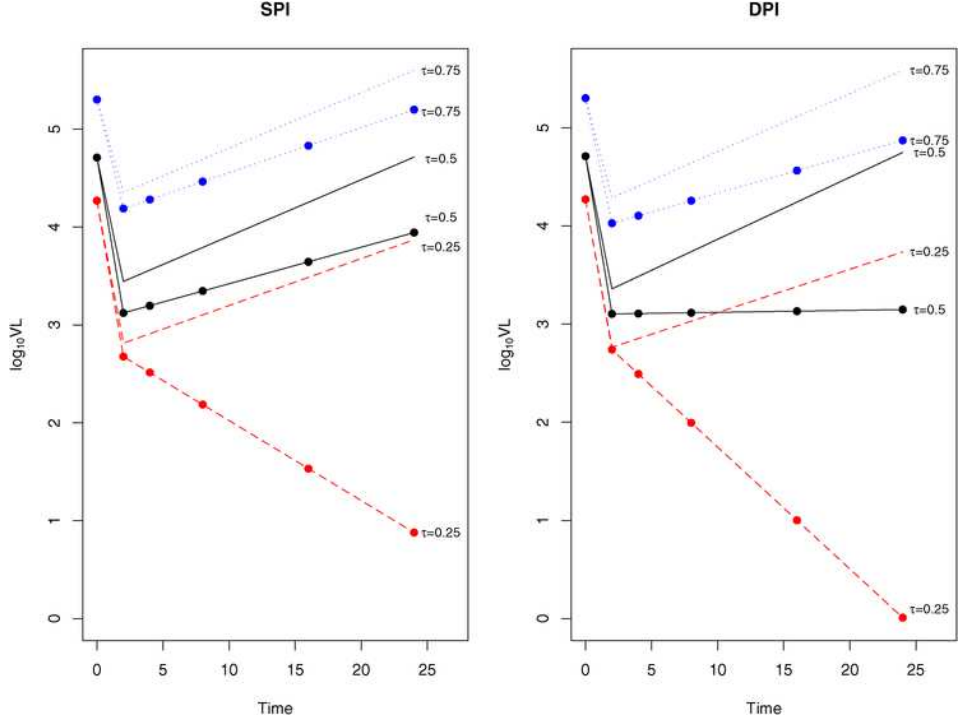

Fig. 3. *Estimated quartiles of* $\log_{10}$ VL *at six visits. The curves with solid dots are from the censored regression and those without dots are from Naïve$_2$ method. The dashed, solid and dotted curves are for* $\tau = 0.25, 0.5$ *and* $0.75$, *respectively.*

To analyze the treatment effect and demonstrate the importance of accounting for censoring and dependency in the data, we tested a series of hypotheses. Table 4 describes these hypotheses and their $p$-values from *QRS*, *Naïve$_2$* and *Indep* at several $\tau$'s. Note that these significance results are for individual $\tau$'s, but not from simultaneous tests. For this data, the $\delta$ is estimated to be 0.37 at median corresponding to a sign correlation of $(0.37 - 0.25)/0.25 = 0.48$. The following are highlights of the interesting findings in the table.

- At $\tau = 0.4$ and $0.5$, our method and the *Indep* method indicate that the treatment is significantly better than the control after week 2. By contrast, the *Naïve$_2$* method indicates no significant difference.

TABLE 4
*Results of inference on $\gamma_1$ and $\gamma_2$ at several quantile levels*

| $\tau$ | Coefficient estimation | | | $p$-value | | | |
|---|---|---|---|---|---|---|---|
| | **Parameter** | **QRS** | **Naïve$_2$** | **Hypothesis** | **QRS** | **Naïve$_2$** | **Indep** |
| 0.4 | $\gamma_1$ | 0.067 | 0.001 | $H_0: \gamma_1 = \gamma_2 = 0$ | 0.022 | 0.965 | 0.100 |
| | $\gamma_2$ | $-0.118$ | 0.002 | $H_0: \gamma_1 = 0$ | 0.126 | 0.930 | 0.066 |
| | | | | $H_0: \gamma_2 = 0$ | 0.001 | 0.841 | 0.015 |
| 0.5 | $\gamma_1$ | $-0.009$ | $-0.042$ | $H_0: \gamma_1 = \gamma_2 = 0$ | 0.001 | 0.717 | 0.001 |
| | $\gamma_2$ | $-0.035$ | 0.005 | $H_0: \gamma_1 = 0$ | 0.900 | 0.507 | 0.872 |
| | | | | $H_0: \gamma_2 = 0$ | 0 | 0.672 | 0.002 |
| 0.6 | $\gamma_1$ | $-0.058$ | 0.033 | $H_0: \gamma_1 = \gamma_2 = 0$ | 0.192 | 0.844 | 0.067 |
| | $\gamma_2$ | $-0.016$ | $-0.009$ | $H_0: \gamma_1 = 0$ | 0.494 | 0.776 | 0.377 |
| | | | | $H_0: \gamma_2 = 0$ | 0.045 | 0.583 | 0.150 |
| 0.7 | $\gamma_1$ | $-0.051$ | $-0.047$ | $H_0: \gamma_1 = \gamma_2 = 0$ | 0.291 | 0.653 | 0.213 |
| | $\gamma_2$ | $-0.009$ | 0.002 | $H_0: \gamma_1 = 0$ | 0.708 | 0.433 | 0.628 |
| | | | | $H_0: \gamma_2 = 0$ | 0.173 | 0.939 | 0.335 |

- At $\tau = 0.6$, our method indicates that the treatment is significantly more favorable (at the 5% level) than the control after week 2. By contrast, neither *Naïve$_2$* nor *Indep* indicates a significant difference.
- For the quantile level $\tau = 0.7$, there is no significant difference between the treatment and the control throughout the trial period—a conclusion supported by all three methods.

Figure 4 highlights the importance of fitting a variety of quantile models to the data. The solid line with open circles in each panel depicts the point estimates, with the shaded area representing a 95% pointwise confidence band following our method. The dashed line represents the mean effect obtained from Tobin's normal censored regression model, with two dotted lines representing a 95% pointwise confidence band for that effect. The confidence interval for the mean effects is computed using the bootstrap method and treating each subject as a single unit. In [28], the $p$-value for testing the interaction effect was 0.0228, suggesting that the $\log_{10}$VL of DPI drops significantly faster than that of SPI throughout the trial period. Our method, however, suggests that two regimens do not differ during the first two weeks across all the quantile levels [Figure 4(a)]. From week 2 to week 24, the mean regression method shows that treatment is more favorable than the control. Quantile regression indicates that this difference mainly comes from the lower tail of the $\log_{10}$VL distribution ($0.25 < \tau \leq 0.45$).

**5. Discussion and conclusions.** In this paper we introduced inference procedures for longitudinal data with fixed censoring within the robust framework of a semi-parametric quantile regression model. One main focus was on providing test procedures for comparing treatments and for assessing the influence of a subset of the covariates. Our proposed quantile rank score test avoids the need for estimating an unknown density. It is relatively easy to implement and performs well in empirical investigations. In particular, by applying our test statistics to data from an AIDS-related clinical trial, we demonstrated the importance of separately considering the various quantiles when assessing the relative merits of different treatments. Moreover, our conclusions could not have been reached by methods that ignore either censoring or intra-subject dependency in the data.

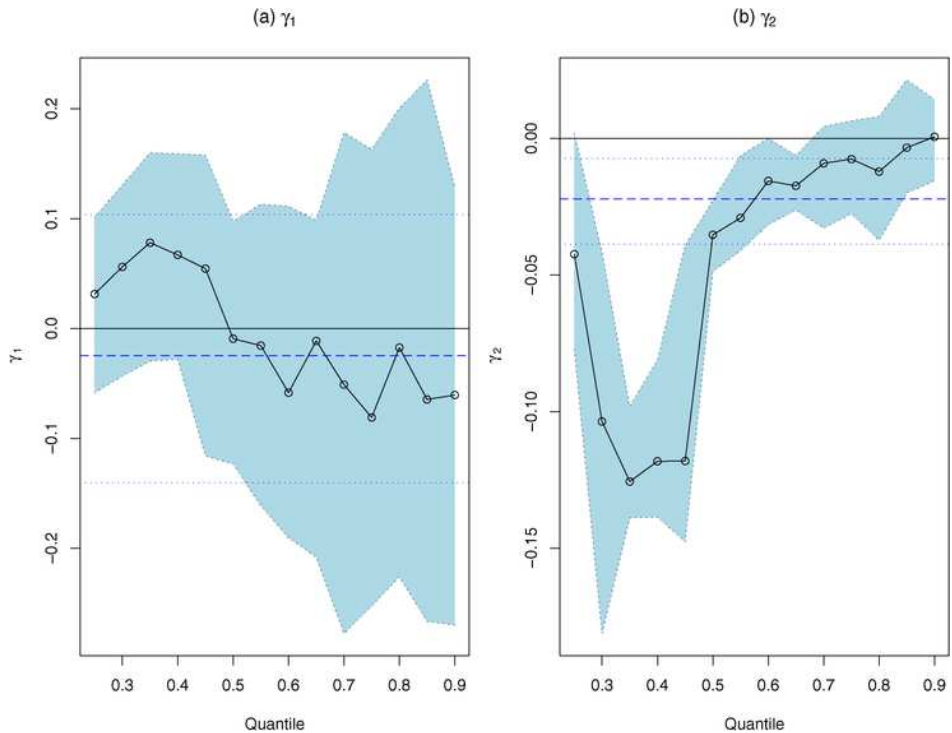

Fig. 4. *The quantile estimates (open circles) of $\gamma_1$ and $\gamma_2$. In each panel, the shaded area depicts a 95% pointwise confidence band for the quantile coefficient, the dashed line represents the mean coefficient estimate with two dotted lines representing a 95% pointwise confidence band for the mean effect.*

The quantile estimate $\hat{\alpha}$ we employed in the current paper is derived under the working assumption of independence. Efficiency might be gained by incorporating appropriate weights to account for the intra-subject correlation structure, as done in [17] for uncensored data. However, He, Fu and Fung [11] and Yin and Cai [33] found, albeit in different contexts, that the efficiency gain in doing so is minimal in finite samples unless the intra-subject correlation is extremely high.

In Section 3, we explored, via Case 4, the robustness of our proposed inference procedure to the assumption of homogeneity of the error terms. Our simulations indicate that the proposed procedures, and the quantile rank test in particular, perform robustly. These findings are encouraging. In a future study, we plan to generalize our methodology to cover models with heteroscedastic errors. In the mean time, it is clear from our proof in the Appendix that for such models our estimator for $\alpha_0$ is still strongly consistent and asymptotically normal, but with a modified asymptotic variance–covariance matrix. Without an appropriate weighting, however, the limiting distribution of our rank score test is no longer that of a chi-squared distribution.

In this paper, the focus was on an exchangeable correlation structure where $\delta = P(u_{ij} < 0, u_{ij'} < 0)$ is taken to be the same for all pairs of errors. Indeed, one can apply the methodology developed here to situations with more general intra-subject dependency structures, as long as one can obtain a consistent estimator for the variance–covariance matrix $\text{Cov}(\varphi(U))$.

## APPENDIX

**A.1. Proof of Theorem 2.1.** We first show the strong consistency of $\hat{\alpha}_0$ by modifying the proof of Theorem 1 in [24] to cover censored regression models for longitudinal data. By the continuity of $Q_n(\alpha)$ defined in (2.6), it suffices to show the strong consistency of $\hat{\alpha}_0$ under $H_0$.

Define
$$u_{ij}^* = y_{ij} - \max(0, x_{ij}^T \alpha_0) = \max(0, x_{ij}^T \alpha_0 + u_{ij}) - \max(0, x_{ij}^T \alpha_0)$$

and
$$h_{ij}(\alpha, \alpha_0) = \max(0, x_{ij}^T \alpha) - \max(0, x_{ij}^T \alpha_0).$$

Note that the minimization of $Q_n(\alpha)$ is equivalent to the minimization of
$$\begin{aligned}
Q_n^*(\alpha) &= n^{-1} \sum_{i=1}^N \sum_{j=1}^{n_i} \rho_\tau \{y_{ij} - \max(0, x_{ij}^T \alpha)\} \\
&\quad - n^{-1} \sum_{i=1}^N \sum_{j=1}^{n_i} \rho_\tau \{y_{ij} - \max(0, x_{ij}^T \alpha_0)\} \\
&= n^{-1} \sum_{i=1}^N \sum_{j=1}^{n_i} \{\rho_\tau(u_{ij}^* - h_{ij}) - \rho_\tau(u_{ij}^*)\} \\
&= n^{-1} \sum_{i=1}^N \sum_{j=1}^{n_i} R_{ij}.
\end{aligned}$$

Then each $R_{ij} = \{\rho_\tau(u_{ij}^* - h_{ij}) - \rho_\tau(u_{ij}^*)\}$ is bounded by $\|x_{ij}\|(\|\alpha\| + \|\alpha_0\|) = O(\|x_{ij}\|)$. By condition A2, Lévy's theorem [26] and Lemma 2.2 of [31], $\hat{\alpha}_0$ will be strongly consistent if the conditional expectation of $Q_n^*(\alpha)$ given the covariates

$$(A.1) \qquad E[Q_n^*(\alpha)|\{x_{ij}\}] \equiv \bar{Q}_n(\alpha)$$

is strictly positive for $\|\alpha - \alpha_0\| \geq \varepsilon$ for arbitrary $\varepsilon > 0$ and all $n$ sufficiently large.

From the derivation of the conditional expectation in (A.1), it can be shown that

$$(A.2) \quad \begin{aligned} \bar{Q}_n(\alpha) &\geq n^{-1} \sum_{ij} I(x_{ij}^T \alpha_0 \geq 0, x_{ij}^T \alpha \geq 0) \int_0^{x_{ij}^T \Delta} (x_{ij}^T \Delta - \lambda) f(\lambda) \, d\lambda \\ &\quad + n^{-1} \sum_{ij} I(x_{ij}^T \alpha_0 \geq 0, x_{ij}^T \alpha < 0) \int_{-x_{ij}^T \alpha_0}^0 (\lambda + x_{ij}^T \alpha_0) f(\lambda) \, d\lambda, \end{aligned}$$

where $\Delta = \alpha - \alpha_0$.

Following the same argument as in (A.10) of Powell [24], condition A4 and inequality (A.2) yield

$$(A.3) \qquad \bar{Q}_n(\alpha) \geq 1/2 \varrho_1 c^2 n^{-1} \sum_{ij} I(x_{ij}^T \alpha_0 \geq \varepsilon_0) I(|x_{ij}^T \Delta| \geq c)$$

for any positive number $c < \min(\varepsilon_0, \varrho_1)$, where $\varepsilon_0$ is defined in condition A3. This completes the proof of strong consistency because the inequality in (A.3) is strictly positive by conditions A2 and A3.

The asymptotic normality follows from the application of Liapunov's central limit theorem to the following lemma:

LEMMA A.1. *Under (2.5) and conditions* A1–A6, *we have*

$$(A.4) \quad \begin{aligned} n^{1/2}(\hat{\alpha}_0 - \alpha_0) &= n^{-1/2} \{f(0)\}^{-1} D_1^{-1} \sum_{ij} I(x_{ij}^T \alpha_0 > 0) x_{ij} \varphi_\tau(u_{ij}) \\ &\quad + D_1^{-1} D_2 \beta_0 + o_p(1). \end{aligned}$$

PROOF. Define

$$\psi_i(\alpha) = \sum_{j=1}^{n_i} I(x_{ij}^T \alpha > 0) x_{ij} \varphi_\tau \{u_{ij} + n^{-1/2} z_{ij}^T \beta_0 - x_{ij}^T (\alpha - \alpha_0)\},$$

$$\Lambda_n(\alpha) = \sum_{i=1}^N E\{\psi_i(\alpha)\}, \qquad u_i(\alpha, d) = \sup_{\|\gamma - \alpha\| \leq d} \|\psi_i(\gamma) - \psi_i(\alpha)\|.$$

Note that $\psi_i(\alpha)$ are independent and

(A.5) $$\psi_i(\hat{\alpha}_0) = \sum_{j=1}^{n_i} I(x_{ij}^T \hat{\alpha}_0 > 0) x_{ij} \varphi_\tau(y_{ij} - x_{ij}^T \hat{\alpha}_0).$$

Since $\hat{\alpha}_0$ minimizes the objective function $Q_n(\alpha)$ in (2.6) for all $\alpha$, the directional derivative of $Q_n(\alpha)$ at $\alpha$ in the direction of a unit vector $\mathbf{v}$ is nonnegative. That is, $\lim_{h \downarrow 0} \{Q_n(\hat{\alpha}_0 + h\mathbf{v}) - Q_n(\hat{\alpha}_0)\}/h \geq 0$. This implies that

$$\sum \sum_{y_{ij} \neq x_{ij}^T \hat{\alpha}_0} I(x_{ij}^T \hat{\alpha}_0 > 0) \varphi_\tau(y_{ij} - x_{ij}^T \hat{\alpha}_0) \mathbf{v}^T x_{ij}$$
$$\leq -\sum \sum_{y_{ij} = x_{ij}^T \hat{\alpha}_0} I(x_{ij}^T \mathbf{v} > 0) \varphi_\tau(y_{ij} - x_{ij}^T \mathbf{v}) \mathbf{v}^T x_{ij}.$$

However, the right-hand side is bounded by

(A.6) $$\sum \sum_{y_{ij} = x_{ij}^T \hat{\alpha}_0} \|x_{ij}\|.$$

By conditions A4 and A5 and the strong consistency of $\hat{\alpha}_0$, we have almost surely a finite number of observations with zero residuals and therefore the quantity in (A.6) is equal to $o(n^{1/4} \gamma_n)$ by condition A2, where $\gamma_n$ is a sequence of positive numbers going to infinity.

Thus, we have

(A.7) $$\sum_{i=1}^{N} \psi_i(\hat{\alpha}_0) = o(n^{1/4} \gamma_n).$$

Since (A.7) is the same as (2.1) in [12], the Bahadur representation in (A.4) will follow from their Theorem 2.2 if conditions B1–B4, B5' and B8 in that theorem hold.

B1. The measurability is trivially satisfied.
B2. This follows directly from the strong consistency of $\hat{\alpha}_0$.
B3. With some manipulations, we obtain

$$u_i(\alpha, d) \leq \sup_{\|\gamma - \alpha\| \leq d} \tau \sum_j \|x_{ij}\| I\{|x_{ij}^T \alpha| \leq |x_{ij}^T(\alpha - \gamma)|\}$$
$$+ \sup_{\|\gamma - \alpha\| \leq d} \sum_j \|x_{ij}\| I\{|u_{ij} - x_{ij}^T \alpha + n^{-1/2} z_{ij}^T \beta_0| \leq |x_{ij}^T(\alpha - \gamma)|\}$$
$$\leq \tau C_{i1} + B_i,$$

where $B_i = \sum_j \|x_{ij}\| I(|u_{ij} - x_{ij}^T\alpha + n^{-1/2}z_{ij}^T\beta_0| \le \|x_{ij}\|d)$ and $C_{i1} = \sum_j \|x_{ij}\|$ $I(|x_{ij}^T\alpha| \le \|x_{ij}\|d)$. Denote $C_{i2} = \sum_j \|x_{ij}\|^2 I\{|x_{ij}^T\alpha| \le \|x_{ij}\|d\}$. Then, by conditions A2 and A4, we have

$$E(B_i) = \sum_j \|x_{ij}\| P\{|u_{ij} - x_{ij}^T\alpha + n^{-1/2}z_{ij}^T\beta_0| \le \|x_{ij}\|d\}$$

$$= \sum_j \|x_{ij}\|^2 \varrho_2 d + O(d^2) = W_{i1}d + O(d^2),$$

$$E(B_i^2) \le n_i \sum_j \|x_{ij}\|^3 \varrho_2 d = W_{i2}d + O(d^2),$$

where $W_{i1}$ and $W_{i2}$ are some positive constants. Therefore, condition B3 follows with $a_i = d^{-1}\tau^2 n_i C_{i2} + 2\tau C_{i1} W_{i1} + W_{i2}$.

B4. Under condition A2, we have $A_n = \sum_i a_i = O(n)$. Thus $A_{2n} = O(A_n)$.

B5′. By condition A2, $\max_i u_i \le n_i \max_{ij} \|x_{ij}\| = O(n^{1/4})$. B5′ follows by taking $d_n = n^{-1/2}(\log n)^4$.

B8. Note that

$$\Lambda_n(\alpha) = \sum_{ij} I(x_{ij}^T\alpha > 0)x_{ij}[\tau - F\{x_{ij}^T(\alpha - \alpha_0) - n^{-1/2}z_{ij}^T\beta_0\}]$$

(A.8)
$$= -nf(0)D_{1n}(\alpha)(\alpha - \alpha_0) + n^{1/2}f(0)D_{2n}(\alpha)\beta_0 + O(1).$$

Therefore we have

(A.9) $$\Lambda_n(\hat{\alpha}_0) = -nf(0)D_{1n}(\hat{\alpha}_0)(\hat{\alpha}_0 - \tilde{\alpha}_0) + O(1),$$

where $\tilde{\alpha}_0 = \alpha_0 + n^{-1/2}\{D_{1n}(\hat{\alpha}_0)\}^{-1}D_{2n}(\hat{\alpha}_0)\beta_0$. Thus B8 holds with $b_n = O(1)$ and $D_n = -nf(0)D_{1n}(\hat{\alpha}_0)$. It then follows from the consistency of $\hat{\alpha}_0$ and Theorem 2.2 of [12] that

$$n^{1/2}(\hat{\alpha}_0 - \tilde{\alpha}_0) = n^{-1/2}\{(f(0)\}^{-1}\{D_{1n}(\hat{\alpha}_0)\}^{-1}\sum_i \psi_i(\tilde{\alpha}_0) + o_p(1)$$

(A.10)
$$= n^{-1/2}\{f(0)\}^{-1}\{D_{1n}(\hat{\alpha}_0)\}^{-1}\sum_{ij} I(x_{ij}^T\tilde{\alpha}_0 > 0)x_{ij}$$

$$\times \varphi_\tau\{u_{ij} + n^{-1/2}z_{ij}^T\beta_0 - x_{ij}^T(\tilde{\alpha}_0 - \alpha_0)\} + o_p(1).$$

By expanding $\varphi_\tau(u + \varepsilon)$ around $u$, we obtain

$$n^{1/2}(\hat{\alpha}_0 - \tilde{\alpha}_0)$$

$$= n^{-1}D_1^{-1}\sum_{ij} I(x_{ij}^T\alpha_0 > 0)x_{ij}(z_{ij}^T\beta_0 - x_{ij}^T D_1^{-1} D_2\beta_0)$$

(A.11)
$$+ n^{-1/2}\{f(0)\}^{-1}D_1^{-1}\sum_{ij} I(x_{ij}^T\alpha_0 > 0)x_{ij}\varphi_\tau(u_{ij}) + o_p(1)$$

$$= n^{-1/2}\{f(0)\}^{-1}D_1^{-1}\sum_{ij} I(x_{ij}^T\alpha_0 > 0)x_{ij}\varphi_\tau(u_{ij}) + o_p(1).$$

This completes the proof of Lemma A.1.   □

**A.2. Proof of Theorem 2.2.**  The proof of Theorem 2.2 relies on the following two lemmas:

LEMMA A.2.  *Define*
$$S_n^* = n^{-1/2}\sum_{ij} z_{ij}^* I(x_{ij}^T\alpha_0 + n^{-1/2}z_{ij}^T\beta_0 > 0)\varphi_\tau(u_{ij}^*),$$
*where* $u_{ij}^* = y_{ij} - \max(0, x_{ij}^T\alpha_0 + n^{-1/2}z_{ij}^T\beta_0)$. *Then under* $H_n$ *and conditions* A1–A7, *we have* $S_n = S_n^* + n^{-1}f(0)\sum_{ij} I(x_{ij}^T\alpha_0 > 0)z_{ij}^* z_{ij}^T \beta_0 + o_p(1)$.

PROOF.  First note that
$$\begin{aligned}S_n^* &= n^{-1/2}\sum_{ij} z_{ij}^* I(x_{ij}^T\alpha_0 + n^{-1/2}z_{ij}^T\beta_0 > 0)\varphi_\tau\\
&\quad\times\{\max(0, u_{ij} + x_{ij}^T\alpha_0 + n^{-1/2}z_{ij}^T\beta_0)\\
&\qquad - \max(0, x_{ij}^T\alpha_0 + n^{-1/2}z_{ij}^T\beta_0)\}\\
&= n^{-1/2}\sum_{ij} z_{ij}^* I(x_{ij}^T\alpha + n^{-1/2}z_{ij}^T\beta_0 > 0)\varphi_\tau(u_{ij}).\end{aligned}$$

For any fixed $t$ such that $\|t\| \leq C$, we define
$$\begin{aligned}R_i(t) &= \sum_{j=1}^{n_i} z_{ij}^*[I(x_{ij}^T\alpha_0 + n^{-1/2}x_{ij}^T t > 0)\\
&\qquad\times\varphi_\tau\{y_{ij} - \max(0, x_{ij}^T\alpha_0 + n^{-1/2}x_{ij}^T t)\}\\
&\qquad - I(x_{ij}^T\alpha_0 + n^{-1/2}z_{ij}^T\beta_0 > 0)\varphi_\tau(u_{ij})]\\
&= \sum_{j=1}^{n_i} z_{ij}^*\{I(x_{ij}^T\alpha_0 + n^{-1/2}x_{ij}^T t > 0)\\
&\qquad\times\varphi_\tau(u_{ij} + n^{-1/2}z_{ij}^T\beta_0 - n^{-1/2}x_{ij}^T t)\\
&\qquad - I(x_{ij}^T\alpha_0 + n^{-1/2}z_{ij}^T\beta_0 > 0)\varphi_\tau(u_{ij})\}.\end{aligned}$$
(A.12)

By condition A4, each coordinate $R_i^{(k)}(t), k = 1,\ldots,q$, satisfies
$$\begin{aligned}\sum_i \mathrm{Var}\{R_i^{(k)}(t)\} &\leq \sum_i n_i \sum_j E\{R_{ij}(t)\}^2\\
&\leq \sum_i n_i \sum_j \|z_{ij}^*\|^2\{\varrho_2|n^{-1/2}x_{ij}^T t - n^{-1/2}z_{ij}^T\beta_0| + I(|x_{ij}^T\alpha_0| \leq d_n\|x_i\|)\},\end{aligned}$$

where $d_n = n^{-1/2}\|x_{ij}\|^{-1} \max(\|x_{ij}\|, \|z_{ij}\|) \cdot \max(\|\beta_0\|, \|t\|) = O(n^{-1/4})$. It then follows from condition A5 that

$$\sum_i \text{Var}\{R_i^{(k)}(t)\} = O(n^{3/4}). \tag{A.13}$$

Under condition A2, it is clear that $\max_i \|R_i(t)\| \leq C n^{1/4}$ for some constant $C$. It follows from the Hoeffding inequality and the chaining argument that

$$\sup_{\|t\|\leq C} \left\|\sum_i R_i(t) - E\left\{\sum_i R_i(t)\right\}\right\| = O_p(n^{3/8}(\log n)^{1/2}). \tag{A.14}$$

Under conditions A2 and A4 and using the orthogonality of $Z^*$ and $\{I(X\alpha_0 > 0) \otimes \underline{1}_n^T\}X$, we obtain that

$$E\left\{\sum_i R_i(t)\right\}$$
$$= \sum_i \sum_j z_{ij}^* \{I(x_{ij}^T \alpha_0 + n^{-1/2} x_{ij}^T t > 0)(\tau - F(n^{-1/2} x_{ij}^T t - n^{-1/2} z_{ij}^T \beta_0))\}$$
$$= \sum_{ij} z_{ij}^* I(x_{ij}^T \alpha_0 + n^{-1/2} x_{ij}^T t > 0) f(0)(n^{-1/2} z_{ij}^T \beta_0 - n^{-1/2} x_{ij}^T t) + O(1)$$
$$= n^{-1/2} f(0) \sum_{ij} z_{ij}^* I(x_{ij}^T \alpha_0 > 0) z_{ij}^T \beta_0 + O(1).$$

This, together with Theorem 2.1 and (A.14), completes the proof. $\square$

LEMMA A.3. *Under $H_n$ and conditions A1–A7, we have $\hat{\delta} \xrightarrow{P} \delta$, as $n \to \infty$.*

PROOF. Recall that $\hat{\delta} = L^{-1} \sum_{ij \neq j'} I(x_{ij}^T \hat{\alpha}_0 > 0, x_{ij'}^T \hat{\alpha}_0 > 0) I(\hat{u}_{ij} < 0, \hat{u}_{ij'} < 0)$. We define $\delta^* = L^{-1} \sum_{ij \neq j'} I(x_{ij}^T \alpha_0 + n^{-1/2} z_{ij}^T \beta_0 > 0, x_{ij'}^T \alpha_0 + n^{-1/2} z_{ij'}^T \beta_0 > 0) I(u_{ij}^* < 0, u_{ij'}^* < 0)$, where $u_{ij}^* = y_{ij} - \max(0, x_{ij}^T \alpha_0 + n^{-1/2} z_{ij}^T \beta_0)$. Using Lemma 4.1 of [12] and the root-$n$ consistency of $\hat{\alpha}_0$, we can establish that $\hat{\delta} - \delta^* = o_p(1)$. Lemma A.3 thus follows by applying the weak law of large numbers. $\square$

PROOF OF THEOREM 2.2. Denote $R_i = \sum_j I(x_{ij}^T \alpha_0 + n^{-1/2} z_{ij}^T \beta_0) z_{ij}^* \varphi_\tau(u_{ij})$. Note that $S_n^* = n^{-1/2} \sum_i R_i$ is the summation of independent entries. It follows from the Lindberg–Feller central limit theorem (CLT) that

$$\{V_n(\delta; \alpha_0)\}^{-1/2} S_n^* \xrightarrow{D} N(\underline{0}_q, I_q). \tag{A.15}$$

The proof of Theorem 2.2 is therefore complete by combining (A.15), Lemmas A.2 and A.3. $\square$

**Acknowledgments.** We thank two anonymous referees and an Associate Editor for their constructive comments and prompt reviews.

## REFERENCES

[1] BUCHINSKY, M. (1994). Changes in the U.S. wage structure 1963–1987: Applications of quantile regression. *Econometrica* **62** 405–458.
[2] BILIAS, Y., CHEN, S. and YING, Z. (2000). Simple resampling methods for censored regression quantiles. *J. Econometrics* **99** 373–386. MR1792254
[3] CHEN, L. and WEI, Y. (2005). Computational issues for quantile regression. *Sankhyā* **67** 399–417. MR2208896
[4] CHU, H., MOULTON, L. H., MACK, W. J., PASSARO, D. J., BARROSO, P. F. and MUÑOZ, A. (2005). Correlating two continuous variables subject to detection limits in the context of mixture distributions. *J. Roy. Statist. Soc. Ser. C* **54** 831–845. MR2209034
[5] DING, A. A. and WU, H. (2001). Assessing antiviral potency of anti-HIV therapies in vivo by comparing viral decay rates in viral dynamic models. *Biostatistics* **2** 13–29.
[6] FITZENBERGER, B. (1997). A guide to censored quantile regressions. In *Handbook of Statistics: Robust Inference* **15** (G. S. Maddala and C. R. Rao, eds.) 405–437. North-Holland, Amsterdam. MR1492720
[7] GHOSH, P., BRANCO, M. D. and CHAKRABORTY, H. (2007). Bivariate random effect model using skew-normal distribution with application to HIV-RNA. *Statistics in Medicine* **26** 1255–1267. MR2345719
[8] GRAY, L., CORTINA-BORJA, M. and NEWELL, M. L. (2004). Modelling HIV-RNA viral load in vertically infected children. *Statistics in Medicine* **23** 769–781.
[9] GUTENBRUNNER, C., JUREĈKOVÁ, J., KOENKER, R. and PORTNOY, S. (1993). Tests of linear hypotheses based on regression rank scores. *J. Nonparametr. Statist.* **2** 307–333. MR1256383
[10] HAMMER, S. M., VAIDA, F., BENNETT, K. K., HOLOHAN, M. K., SHEINER, L., ERON, J. J., WHEAT, L. J., MITSUYASU, R. T., GULICK, R. M., VALENTINE, F. T., ABERG, J. A., ROGERS, M. D., KAROL, C. N., SAAH, A. J., LEWIS, R. H., BESSEN, L. J., BROSGART, C., DEGRUTTOLA, V. and MELLORS, J. W. (2002). Dual vs single protease inhibitor therapy following antiretroviral treatment failure: A randomized trial. *J. Amer. Med. Assoc.* **288** 169–180.
[11] HE, X., FU, B. and FUNG, W. K. (2003). Median regression for longitudinal data. *Statistics in Medicine* **22** 3655–3669.
[12] HE, X. and SHAO, Q. M. (1996). A general Bahadur representation of $M$-estimators and its application to linear regression with nonstochatic designs. *Ann. Statist.* **24** 2608–2630. MR1425971
[13] HOGAN, J. W. and LEE, J. Y. (2004). Marginal structural quantile models for longitudinal observational studies with time-varying. *Statist. Sinica* **14** 927–944. MR2089340
[14] HUANG, W., DE GRUTTOLA, V., FISCHL, M., HAMMER, S., RICHMAN, D., HAVLIR, D., GULICK, R., SQUIRES, K. and MELLORS, J. (2001). Patterns of plasma human immunodeficiency virus type 1 RNA response to antiretroviral therapy. *J. Infectious Diseases* **183** 1455–1465.
[15] HUGHES, J. P. (1999). Links mixed effects models with censored data with application to HIV RNA levels. *Biometrics* **55** 625–629.

[16] JACQMIN-GADDA, H., THIÉBAUT, R., CHÊNE, G. and COMMENGES, D. (2000). Analysis of left censored longitudinal data with application to viral load in HIV infection. *Biostatistics* **1** 355–368.

[17] JUNG, S. (1996). Quasi-likelihood for median regression models. *J. Amer. Statist. Assoc.* **91** 251–257. MR1394079

[18] KOCHERGINSKY, M., HE, X. and MU, Y. (2005). Practical confidence intervals for regression quantiles. *J. Comput. Graph. Statist.* **14** 41–55. MR2137889

[19] KOENKER, R. (1994). Confidence intervals for regression quantiles. In *Asymptotic Statistics: Proceedings of the 5th Prague Symposium on Asymptotic Statistics* (P. Mandl and M. Husková, eds.) 349–359. Physica-Verlag, Heidelberg. MR1311953

[20] KOENKER, R. and PARK, B. (1996). An interior point algorithm for nonlinear quantile regression. *J. Econometrics* **71** 265–283. MR1381085

[21] LYLES, R. H., LYLES, C. M. and TAYLOR, D. J. (2000). Random regression models for human immunodeficiency virus ribonucleic acid data subject to left censoring and informative drop-outs. *J. Roy. Statist. Soc. Ser. C* **49** 485–497. MR1824554

[22] MOULTON, L. H. and HALSEY, N. A. (1995). A mixture model with detection limits for regression analyses of antibody response to vaccine. *Biometrics* **51** 1570–1578.

[23] PORTNOY, S. (2003). Censored regression quantiles. *J. Amer. Statist. Assoc.* **98** 1001–1012. MR2041488

[24] POWELL, J. L. (1984). Least absolute deviations estimation for the censored regression model. *J. Econometrics* **25** 303–325. MR0752444

[25] POWELL, J. L. (1986). Censored regression quantiles. *J. Econometrics* **32** 143–155. MR0853049

[26] RESNICK, S. (2003). *A Probability Path*. Birkhäuser, Boston. MR1664717

[27] SINGH, A. and NOCERINO, J. (2002). Robust estimation of mean and variance using environmental data sets with below detection limit observations. *Chemometrics Intelligent Laboratory Systems* **60** 69–86.

[28] SUN, Y. and WU, H. (2005). Semiparametric time-varying coefficients regression model for longitudinal data. *Scand. J. Statist.* **32** 21–47. MR2136800

[29] WANG, H. and HE, X. (2007). Detecting differential expressions in GeneChip microarray studies: A quantile approach. *J. Amer. Statist. Assoc.* **102** 104–112. MR2293303

[30] WEI, Y. and HE, X. (2006). Conditional growth charts. *Ann. Statist.* **34** 2069–2097. MR2291494

[31] WHITE, H. (1980). Nonlinear regression on cross-sectional data. *Econometrica* **48** 721–746. MR0573321

[32] WU, L. (2002). A joint model for nonlinear mixed-effects models with censoring and covariates measured with error, with application to AIDS studies. *J. Amer. Statist. Assoc.* **97** 955–964. MR1951254

[33] YIN, G. S. and CAI, J. (2005). Quantile regression models with multivariate failure time data. *Biometrics* **61** 151–161. MR2135855

[34] YING, Z., JUNG, S. H. and WEI, L. J. (1995). Survival analysis with median regression models. *J. Amer. Statist. Assoc.* **90** 178–184. MR1325125

[35] ZHAO, L. C. (2004). Linear hypothesis testing in censored regression models. *Statist. Sinica* **14** 333–347. MR2036776

DEPARTMENT OF STATISTICS  
NORTH CAROLINA STATE UNIVERSITY  
RALEIGH, NORTH CAROLINA 27695  
USA  
E-MAIL: wang@stat.ncsu.edu

INFORMATION AND OPERATIONS MANAGEMENT  
UNIVERSITY OF SOUTHERN CALIFORNIA  
LOS ANGELES, CALIFORNIA 90089  
USA  
E-MAIL: fygenson@marshall.usc.edu